# Zeta-functions for germs of meromorphic functions and Newton diagrams

S.M. Gusein-Zade    I. Luengo    A. Melle-Hernández [*]

**Abstract**

For a germ of a meromorphic function $f = \frac{P}{Q}$, we offer notions of the monodromy operators at zero and at infinity. If the holomorphic functions $P$ and $Q$ are non-degenerated with respect to their Newton diagrams, we give an analogue of the formula of Varchenko for the zeta-functions of these monodromy operators.

## 1 Germs of meromorphic functions

A polynomial $f$ of $(n+1)$ complex variables of degree $d$ determines a meromorphic function $f$ on $\mathbb{CP}^{n+1}$. If one wants to understand the behaviour of $f$ at infinity, it is natural to analize germs of the meromorphic function $f$ at points from the infinite hyperplane $\mathbb{CP}^n_\infty \subset \mathbb{CP}^{n+1}$. In local analytic coordinates $z_0, z_1, \ldots, z_n$, centred at a point $p \in \mathbb{CP}^n_\infty$ such that the infinite hyperplane $\mathbb{CP}^n_\infty$ is given by the equation $\{z_0 = 0\}$, the germ of the function $f$ has the form $f = \dfrac{P(z_0, \ldots, z_n)}{z_0^d}$. Let us consider germs of meromorphic functions of a general form.

DEFINITION 1 A *germ of a meromorphic function on* $(\mathbb{C}^{n+1}, 0)$ is a fraction $f = \frac{P}{Q}$, where $P$ and $Q$ are germs of holomorphic functions $(\mathbb{C}^{n+1}, 0) \to (\mathbb{C}, 0)$. Two germs of meromorphic functions $f = \frac{P}{Q}$ and $f' = \frac{P'}{Q'}$ are *equal* if there exists a germ of a holomorphic function $U : (\mathbb{C}^{n+1}, 0) \to \mathbb{C}$ such that $U(0) \neq 0$, $P' = U \cdot P$ and $Q' = U \cdot Q$.

REMARKS. (1) For our convenience here we do not consider functions of the type $\frac{1}{Q(z)}$ or $\frac{P(z)}{1}$.

(2) According to the definition $\dfrac{x}{y} \neq \dfrac{x^2}{xy}$, but $\dfrac{x}{y} = \dfrac{x \exp(x)}{y \exp(x)}$.

Recently V.I. Arnold had obtained classifications of simple germs of meromorphic functions for certain equivalence relations.

In what follows we shall consistently use resolutions of germs of meromorphic functions.

[*]First author is partially supported by RFBR-95-01-01122a and INTAS-4373. Last two authors are partially supported by CAICYT PB94-291.

DEFINITION 2 A *resolution of the germ* $f$ is a modification of $(\mathbb{C}^{n+1}, 0)$ (i.e. a proper analytic map $\pi : \mathcal{X} \to \mathcal{U}$ of a smooth analytic manifold $\mathcal{X}$ onto a neighbourhood $\mathcal{U}$ of the origin in $\mathbb{C}^{n+1}$, which is an isomorphism outside of a proper analytic subspace in $\mathcal{U}$) such that the total transform $\pi^{-1}(H)$ of the hypersurface $H = \{P = 0\} \cup \{Q = 0\}$ is a normal crossing divisor at each point of $\mathcal{X}$.

The fact that the preimage $\pi^{-1}(H)$ is a divisor with normal crossings implies that in a neighbourhood of any point of it, there exists a local system of coordinates $y_0, y_1, \ldots, y_n$ such that the liftings $\tilde{P} = P \circ \pi$ and $\tilde{Q} = Q \circ \pi$ of the functions $P$ and $Q$ to the space $\mathcal{X}$ of the modification are equal to $u\, y_0^{k_0} y_1^{k_1} \cdots y_n^{k_n}$ and $v\, y_0^{l_0} y_1^{l_1} \cdots y_n^{l_n}$ respectively, where $u(0) \neq 0$ and $v(0) \neq 0$.

Let $B_\varepsilon$ be the closed ball of radius $\varepsilon$ with the centre at the origin in $\mathbb{C}^{n+1}$ and $\varepsilon$ be small enough such that (representatives of) the functions $P$ and $Q$ are defined in $B_\varepsilon$ and for any positive $\varepsilon' < \varepsilon$ the sphere $S_{\varepsilon'} = \partial B_{\varepsilon'}$ intersects the analytic spaces $\{P = 0\}$, $\{Q = 0\}$ and $\{P = Q = 0\}$ transversally (in the stratified sense). We choose $\delta$ small enough and take the ball $B_\delta \subset \mathbb{C}^2$ of radius $\delta$ centred at the origin.

DEFINITION 3 Let $c \in \mathbb{C}$ be such that $\|c\|$ is small enough, the 0-*Milnor fibre* $\mathcal{M}_f^0$ of the germ $f$ is the set

$$\mathcal{M}_f^0 = \{z \in B_\varepsilon : (P(z), Q(z)) \in B_\delta \subset \mathbb{C}^2,\ f(z) = \frac{P(z)}{Q(z)} = c\}.$$

In the same way, for $c \in \mathbb{C}$ such that $\|c\|$ is large enough, the $\infty$-*Milnor fibre* $\mathcal{M}_f^\infty$ of the germ $f$ is the set

$$\mathcal{M}_f^\infty = \{z \in B_\varepsilon : (P(z), Q(z)) \in B_\delta \subset \mathbb{C}^2,\ f(z) = \frac{P(z)}{Q(z)} = c\}.$$

LEMMA 1 *The notion of the* 0- *(respectively of the* $\infty$-*) Milnor fibre is well defined, i.e. for* $\|c\|$ *small enough:* $\|c\| \ll \delta \ll \varepsilon$ *(respectively for* $\|c\|$ *large enough:* $\|c\|^{-1} \ll \delta \ll \varepsilon$*) the differentiable type of* $\mathcal{M}_f^0$ *(respectively of* $\mathcal{M}_f^\infty$*) does not depend on* $\varepsilon, \delta$ *and* $c$.

*Proof.* Let $\pi : \mathcal{X} \to \mathcal{U}$ be a resolution of the germ $f$ which is an isomorphism outside the hypersurface $H = \{P = 0\} \cup \{Q = 0\}$. Let $r : \mathbb{C}^{n+1} \to \mathbb{R}$ be the function $r(z) = \|z\|^2$, let $\tilde{r} = r \circ \pi : \mathcal{X} \to \mathbb{R}$ be the lifting of $r$ to the space $\mathcal{X}$ of the resolution. For $\varepsilon$ small enough, the hypersurface $\tilde{S}_\varepsilon = \{\tilde{r} = \varepsilon^2\}$ (the preimage of the sphere $S_\varepsilon \subset \mathbb{C}^{n+1}$) is transversal to all components of the total transform $\pi^{-1}(H)$. At each point of $\pi^{-1}(H)$ in a local coordinate system one has $P \circ \pi = u y_0^{k_0} \cdots y_n^{k_n}$, $Q \circ \pi = v y_0^{l_0} \cdots y_n^{l_n}$ with $u(0) \neq 0$ and $v(0) \neq 0$. Thus $f \circ \pi = w y_0^{m_0} \cdots y_n^{m_n}$ with $w(0) \neq 0$. The real hypersurface $\tilde{S}_\varepsilon$ is transversal to all coordinate subspaces (of different dimensions). It is not difficult to show that this implies transversality of $\tilde{S}_\varepsilon$ to the (complex) hypersurfaces $\{w y_0^{m_0} \cdots y_n^{m_n} = c\}$ for $\|c\|$ small enough and for $\|c\|$ large enough. Now the proof follows from the standard arguments.

REMARKS. (1) The definition means that $\mathcal{M}_f^0$ or $\mathcal{M}_f^\infty$ is equal to

$$\{z \in B_\varepsilon : (P(z), Q(z)) \in B_\delta \subset \mathbb{C}^2, P(z) = cQ(z), \quad P(z) \neq 0\}$$

and thus the Milnor fibres of the functions $\frac{P}{Q}$ and $\frac{RP}{RQ}$ with $R(0) = 0$ are, generally speaking, different.

(2) For $f = \frac{P}{Q}$, let $f^{-1} = \frac{Q}{P}$. It is not difficult to understand that $\mathcal{M}_{f^{-1}}^0 = \mathcal{M}_f^\infty$ and $\mathcal{M}_{f^{-1}}^\infty = \mathcal{M}_f^0$. Just the same properties hold for the monodromy transformations and for the zeta-functions discussed bellow.

(3) It is possible (and sometimes more convenient) to define the Milnor fibres as follows:

$$\mathcal{M}_f^0 = \{z \in B_\varepsilon : \|Q(z)\| \leq \delta, P(z) = cQ(z) \neq 0\}$$

with $\|c\| \ll \delta \ll \varepsilon$, and

$$\mathcal{M}_f^\infty = \{z \in B_\varepsilon : \|P(z)\| \leq \delta, P(z) = cQ(z) \neq 0\}$$

with $\|c\|^{-1} \ll \delta \ll \varepsilon$.

The meromorphic function $f$ determines a map from $B_\varepsilon \setminus \{P = Q = 0\}$ to the projective line $\mathbb{CP}^1$ ($z \mapsto (P(z) : Q(z))$), which also will be denoted by $f$. Lemma 1 implies that this map is a locally trivial fibration in punctured neighbourhoods of the points $0 = (0 : 1)$ and $\infty = (1 : 0)$ of $\mathbb{CP}^1$.

DEFINITION 4 The 0-*monodromy transformation* $h_f^0$ (respectively the $\infty$-*monodromy transformation* $h_f^\infty$) of the germ $f$ is the monodromy transformation of the fibration $f$ over the loop $c \cdot \exp(2\pi it)$, $t \in [0, 1]$, with $\|c\|$ small enough (respectively large enough).

The 0- or $\infty$- monodromy operator is the action of the corresponding monodromy transformation in a homology group of the Milnor fibre. We are interested to apply the results for meromorphic functions to the problem of calculating the zeta-function of a polynomial at infinity. Thus we shall consider the zeta-functions $\zeta_f^0(t)$ and $\zeta_f^\infty(t)$ of the corresponding monodromy transformations:

$$\zeta_f^\bullet = \prod_{q \geq 0} \{\det [id - t\, h_{f*}^\bullet|_{H_q(\mathcal{M}_f^\bullet; \mathbb{C})}]\}^{(-1)^q}$$

($\bullet = 0$ or $\infty$). This definition coincides with that in [2] and differs by minus sign in the exponent from that in [1].

## 2 Resolution of singularities and the formula of A'Campo for germs of meromorphic functions

Let $f = \frac{P}{Q}$ be a germ of a meromorphic function on $(\mathbb{C}^{n+1}, 0)$ and let $\pi : \mathcal{X} \to \mathcal{U}$ be a resolution of the germ $f$. The preimage $\mathcal{D} = \pi^{-1}(0)$ of the origin of $\mathbb{C}^{n+1}$, is a normal crossing divisor. Let $S_{k,l}$ be the set of points of $\mathcal{D}$ in a neighbourhood of which the functions $P \circ \pi$ and $Q \circ \pi$ in some local coordinates have the form $u\, y_0^k$ and $v\, y_0^l$ respectively ($u(0) \ne 0$, $v(0) \ne 0$). A slight modification of the arguments of A'Campo ([1]) permits to obtain the following version of his formula for the zeta-function of the monodromy of a meromorphic function.

THEOREM 1 *Let the resolution $\pi : \mathcal{X} \to \mathcal{U}$ be an isomorphism outside the hypersurface $H = \{P = 0\} \cup \{Q = 0\}$. Then*

$$\zeta_f^0(t) = \prod_{k>l}(1 - t^{k-l})^{\chi(S_{k,l})},$$
$$\zeta_f^\infty(t) = \prod_{k<l}(1 - t^{l-k})^{\chi(S_{k,l})}.$$

REMARK. A resolution $\pi$ of the germ $f' = \frac{RP}{RQ}$ is at the same time a resolution of the germ $f = \frac{P}{Q}$. Moreover the multiplicities of any component $C$ of the exceptional divisor in the zero divisors of the liftings $(RP) \circ \pi$ and $(RQ) \circ \pi$ of the germs $RP$ and $RQ$ are obtained from those of the germs $P$ and $Q$ by adding one and the same integer, the multiplicity $m = m(C)$ of $R$. Nevertheless the meromorphic functions $f$ and $f'$ can have different zeta-functions. The reason why formulae in the previous theorem give different results for $f$ and $f'$ consists in the fact that if an open part of the component $C$ lies in $S_{k,l}(f)$ then, generally speaking, its part which lies in $S_{k+m,l+m}(f')$ is smaller.

## 3 Zeta-functions of meromorphic functiuons via partial resolutions

Let $f = \frac{P}{Q}$ be a germ of a meromorphic function on $(\mathbb{C}^{n+1}, 0)$ and let $\pi : (\mathcal{X}, \mathcal{D}) \to (\mathbb{C}^{n+1}, 0)$ be an arbitrary modification of $(\mathbb{C}^{n+1}, 0)$, which is an isomorphism outside the hypersurface $H = \{P = 0\} \cup \{Q = 0\}$ (i.e. $\pi$ is not necessarily a resolution). Let $\varphi = f \circ \pi$ be the lifting of $f$ to the space of the modification, i.e. the meromorphic function $\frac{P \circ \pi}{Q \circ \pi}$. For a point $x \in \pi^{-1}(H)$, let $\zeta_{\varphi,x}^0(t)$ and $\zeta_{\varphi,x}^\infty(t)$ be the zeta-functions of the 0- and $\infty$-monodromies of the germ of the function $\varphi$ at $x$. Let $\mathcal{S} = \{\Xi\}$ be a prestratification of $\mathcal{D} = \pi^{-1}(0)$ (that is a partitioning into semi-analytic subspaces without any regularity conditions) such that, for each stratum $\Xi$ of $\mathcal{S}$, the zeta-functions $\zeta_{\varphi,x}^0(t)$ and $\zeta_{\varphi,x}^\infty(t)$ do not depend on $x$, for $x \in \Xi$. We denote this zeta-functions by $\zeta_\Xi^0$ and by $\zeta_\Xi^\infty$ respectively. The same arguments which were used in [4] imply

THEOREM 2 *For $\bullet = 0$ or $\infty$,*

$$\zeta_f^\bullet(t) = \prod_{\Xi \in \mathcal{S}}[\zeta_\Xi^\bullet(t)]^{\chi(\Xi)}.$$

# 4 Zeta-functions via Newton diagrams

For a germ $R = \sum a_k x^k : (\mathbb{C}^{n+1}, 0) \to (\mathbb{C}, 0)$ of a holomorphic function ($k = (k_0, k_1, \ldots, k_n)$, $x^k = x_0^{k_0} x_1^{k_1} \cdots x_n^{k_n}$), its Newton diagram $\Gamma = \Gamma(R)$ is the union of the compact faces of the polytope $\Gamma_+ = \Gamma_+(R) = $ convex hull of $\bigcup_{k: a_k \neq 0} (k + \mathbb{R}_+^{n+1}) \subset \mathbb{R}_+^{n+1}$.

Let $f = \frac{P}{Q}$ be a germ of a meromorphic function on $(\mathbb{C}^{n+1}, 0)$ and let $\Gamma_1 = \Gamma(P)$ and $\Gamma_2 = \Gamma(Q)$ be the Newton diagrams of $P$ and $Q$. We call the pair $\Lambda = (\Gamma_1, \Gamma_2)$ of Newton diagrams $\Gamma_1$ and $\Gamma_2$ the *Newton pair of $f$*. We say that the germ of the meromorphic function $f$ is *non-degenerated with respect to its Newton pair* $\Lambda = (\Gamma_1, \Gamma_2)$ if the pair of germs $(P, Q)$ is non-degenerated with respect to the pair $\Lambda = (\Gamma_1, \Gamma_2)$ in the sense of [7] (which is an adaptation for *germs* of complete intersections of the definition of A.G. Khovanskii, [5]).

Let us define zeta-functions $\zeta_\Lambda^0(t)$ and $\zeta_\Lambda^\infty(t)$ for a Newton pair $\Lambda = (\Gamma_1, \Gamma_2)$. Let $1 \leq l \leq n+1$ and let $\mathcal{I}$ be a subset of $\{0, 1, \ldots, n\}$ with the number of elements $\#\mathcal{I}$ equal to $l$. Let $L_\mathcal{I}$ be the coordinate subspace $L_\mathcal{I} = \{k \in \mathbb{R}^{n+1} : k_i = 0 \text{ for } i \notin \mathcal{I}\}$ and $\Gamma_{i,\mathcal{I}} = \Gamma_i \cap L_\mathcal{I} \subset L_\mathcal{I}$. Let $L_\mathcal{I}^*$ be the dual of $L_\mathcal{I}$ and $L_{\mathcal{I}+}^*$ the positive orthant of it ( the set of covectors which have positive values on $L_{\mathcal{I} \geq 0} = \{k \in L_\mathcal{I} : k_i \geq 0 \text{ for } i \in \mathcal{I}\}$). For a primitive integer covector $a \in (\mathbb{R}^*)_+^{n+1}$, let $m(a, \Gamma) = \min_{x \in \Gamma}(a, x)$ and let $\Delta(a, \Gamma) = \{x \in \Gamma : (a, x) = m(a, \Gamma)\}$. We denote by $m_\mathcal{I}$ and $\Delta_\mathcal{I}$ the corresponding objects for the diagram $\Gamma_\mathcal{I}$ and a primitive integer covector $a \in L_{\mathcal{I}+}^*$. Let $E_\mathcal{I}$ be the set of primitive integer covectors $a \in L_{\mathcal{I}+}^*$ such that $\dim(\Delta(a, \Gamma_1) + \Delta(a, \Gamma_2)) = l - 1$ (the Minkowski sum $\Delta_1 + \Delta_2$ of two polytopes $\Delta_1$ and $\Delta_2$ is the polytope $\{x = x_1 + x_2 : x_1 \in \Delta_1, \ x_2 \in \Delta_2\}$). There exists only a finite number of such covectors. For $a \in E_\mathcal{I}$, let $\Delta_1 = \Delta(a, \Gamma_1)$, $\Delta_2 = \Delta(a, \Gamma_2)$ and

$$V_a = \sum_{s=0}^{l-1} V_{l-1}(\underbrace{\Delta_1, \ldots, \Delta_1}_{s \text{ terms}}, \underbrace{\Delta_2, \ldots, \Delta_2}_{l-1-s \text{ terms}}),$$

where the definition of the (Minkowski) mixed volume $V(\Delta_1, \ldots, \Delta_m)$ can be found e.g. in [3] or [7]; $(l-1)$-dimensional volume in a rational $(l-1)$-dimensional affine subspace of $L_\mathcal{I}$ has to be normalized in such a way that the volume of the unit cube spanned by any integer basis of the corresponding linear subspace is equal to 1. Let us recall that $V_m(\underbrace{\Delta, \ldots, \Delta}_{m \text{ terms}})$ is simply the $m$-dimensional volume of $\Delta$. We have to assume that $V_0(\text{nothing}) = 1$, (this is necessary to define $V_a$ for $l = 1$). Let:

$$\zeta_\mathcal{I}^0(t) = \prod_{a \in E_\mathcal{I} : m(a,\Gamma_1) > m(a,\Gamma_2)} (1 - t^{m(a,\Gamma_1) - m(a,\Gamma_2)})^{(l-1)! V_a},$$

$$\zeta_\mathcal{I}^\infty(t) = \prod_{a \in E_\mathcal{I} : m(a,\Gamma_1) < m(a,\Gamma_2)} (1 - t^{m(a,\Gamma_2) - m(a,\Gamma_1)})^{(l-1)! V_a},$$

$$\zeta_l^\bullet(t) = \prod_{\mathcal{I}: \#(\mathcal{I})=l} \zeta_\mathcal{I}^\bullet(t),$$

$$\zeta_\Lambda^\bullet(t) = \prod_{l=1}^{n+1} (\zeta_l^\bullet(t))^{(-1)^{l-1}},$$

where $\bullet = 0$ or $\infty$.

**THEOREM 3** *Let $f = \frac{P}{Q}$ be a germ of a meromorphic function on $(\mathbb{C}^{n+1}, 0)$ non-degenerated with respect to its Newton pair $\Lambda = (\Gamma_1, \Gamma_2)$. Then*

$$\zeta_f^0(t) = \zeta_\Lambda^0(t) \quad and \quad \zeta_f^\infty(t) = \zeta_\Lambda^\infty(t).$$

*Proof.* Let $\Sigma$ be an unimodular simplicial subdivision of $\mathbb{R}_{\geq 0}^{n+1}$ which corresponds to the pair $(\Gamma_1, \Gamma_2)$ of Newton diagrams in the sense of [7] Section 4. This subdivision is consistent with each of the Newton diagrams $\Gamma_1$ and $\Gamma_2$ in the sense of [8].

Let $\pi : (\mathcal{X}, \mathcal{D}) \to (\mathbb{C}^{n+1}, 0)$ be the toroidal modification map corresponding to $\Sigma$. Since the pair $(P, Q)$ is non-degenerated with respect to $(\Gamma_1, \Gamma_2)$, $\pi$ is a resolution of the germ $f = \frac{P}{Q}$ (see [7]). We have the sets $S_{k,l} = S_k(P) \cap S_l(Q)$. The description of $S_k(P)$ (and of $S_l(Q)$) can be found in [8], Section 7. Each of them consists of open parts of certain complex tori of some dimensions.

Tori of dimension $n$ correspond to one-dimensional cone of $\Sigma$ which are positive (i.e., lie in $(\mathbb{R}^*)_+^{n+1}$). The multiplicity of $P \circ \pi$ (respectively of $Q \circ \pi$) at such a torus is equal to $m(a, \Gamma_1)$, (respectively to $m(a, \Gamma_2)$) for the primitive integer covector $a$ which spans the corresponding cone.

Tori of dimension $(l-1)$ correspond to positive simplicial $(n+2-l)$-dimensional cones of $\Sigma$ which have a cone of the form

$$\mathfrak{S} = \{a \in (\mathbb{R}^*)_{\geq 0}^{n+1} : a_j > 0 \text{ for } j \notin \mathcal{I}, \quad a_j = 0 \text{ for } j \in \mathcal{I}\}$$

with $\#(\mathcal{I}) = l$ (these cones are elements of $\Sigma$) as its face. Moreover these cones correspond to one-dimensional cones of a partitioning of $L_\mathcal{I}$ which is consistent with the Newton diagrams $\Gamma_{i,\mathcal{I}} = \Gamma_i \cap L_\mathcal{I} \subset L_\mathcal{I}$. The multiplicities of $P \circ \pi$ and $Q \circ \pi$ at such a torus again are equal to $m_\mathcal{I}(a, \Gamma_{1,\mathcal{I}})$ and $m_\mathcal{I}(a, \Gamma_{2,\mathcal{I}})$ for the primitive integer covector $a$ from the corresponding one-dimensional cone.

In order to apply Theorem 1 we have to calculate the Euler characteristic of the corresponding part of an $l-1$-dimensional torus $T$: the complement to the intersection with the strict transform of the hypersurface $H = \{P = 0\} \cup \{Q = 0\}$. Let $A$ (respectively $B$) the intersection of the torus $T$ with the strict transform of the hypersurface $\{P = 0\}$ (respectively of $\{Q = 0\}$), let $\Delta_i := \Delta(a, \Gamma_{i,\mathcal{I}})$. From the results of Khovanskii ([6]) it follows that the Euler characteristic of $A$ (respectively of $B$) is equal to $(-1)^l(l-1)! V_{l-1}(\underbrace{\Delta_1, \ldots, \Delta_1}_{l-1 \text{ terms}})$ (respectively to $(-1)^l(l-1)! V_{l-1}(\underbrace{\Delta_2, \ldots, \Delta_2}_{l-1 \text{ terms}})$), the Euler charecteristic of $A \cap B$ is equal to

$$(-1)^{l-1}(l-1)! [V_{l-1}(\underbrace{\Delta_1, \ldots, \Delta_1}_{l-2 \text{ terms}}, \Delta_2) + V_{l-1}(\underbrace{\Delta_1, \ldots, \Delta_1}_{l-3 \text{ terms}}, \Delta_2, \Delta_2) + \ldots + V_{l-1}(\Delta_1, \underbrace{\Delta_2, \ldots, \Delta_2}_{l-2 \text{ terms}})].$$

Thus the Euler characteristic of the complement of $A \cup B$ in the torus $T$ is equal to

$$\chi(T) - \chi(A) - \chi(B) + \chi(A \cap B) =$$

$$= (-1)^{l-1}(l-1)![V_{l-1}(\underbrace{\Delta_1, \ldots, \Delta_1}_{l-1 \text{ terms}}) + V_{l-1}(\underbrace{\Delta_1, \ldots, \Delta_1}_{l-2 \text{ terms}}, \Delta_2) + \ldots + V_{l-1}(\underbrace{\Delta_2, \ldots, \Delta_2}_{l-1 \text{ terms}})],$$

which implies the statement. □

# 5 The Varchenko type formula for $f = \frac{P}{z_0^d}$

As we have mentioned at the beginning, in order to study the behaviour of polynomials at infinity, germs of meromorphic functions of the form $\frac{P(z_0, z_1, \ldots, z_n)}{z_0^d}$ have to be of interest. In this case the formulae for the zeta-functions $\zeta_\Lambda^0(t)$ and $\zeta_\Lambda^\infty(t)$ are considerably reduced. Thus let us reformulate the definition of these zeta-functions for the case when the Newton diagram $\Gamma_2$ consists of one point $(d, 0, \ldots, 0)$ (in terms of the Newton diagram $\Gamma := \Gamma_1$ of $P$). The description is as follows.

Let $1 \leq l \leq n+1$ and let $\mathcal{I}$ be a subset of $\{1, \ldots, n\}$ with the number of elements $\#\mathcal{I}$ equal to $l-1$. Let $\gamma_1^\mathcal{I}, \ldots, \gamma_{j(\mathcal{I})}^\mathcal{I}$ be all $(l-1)$-dimensional faces of $\Gamma_{\mathcal{I} \cup \{0\}}$ and $a_{\mathcal{I},1}, \ldots, a_{\mathcal{I},j(\mathcal{I})}$ the corresponding primitive covectors (normal to $\gamma_1^\mathcal{I}, \ldots, \gamma_{j(\mathcal{I})}^\mathcal{I}$), $a_{\mathcal{I},s}^0$ is the 0th coordinate of $a_{\mathcal{I},s}$, let $m_s(\mathcal{I}) = (a_{\mathcal{I},s}, k)$ for $k \in \gamma_s^\mathcal{I}$. Then

$$\zeta_{\mathcal{I} \cup \{0\}}^0(t) = \prod_{1 \leq s \leq j(\mathcal{I}): m_s(\mathcal{I}) > d \cdot a_{\mathcal{I},s}^0} (1 - t^{m_s(\mathcal{I}) - d \cdot a_{\mathcal{I},s}^0})^{(l-1)!V_{l-1}(\gamma_s^\mathcal{I})},$$

$$\zeta_{\mathcal{I} \cup \{0\}}^\infty(t) = \prod_{1 \leq s \leq j(\mathcal{I}): m_s(\mathcal{I}) < d \cdot a_{\mathcal{I},s}^0} (1 - t^{d \cdot a_{\mathcal{I},s}^0 - m_s(\mathcal{I})})^{(l-1)!V_{l-1}(\gamma_s^\mathcal{I})},$$

$$\zeta_l^\bullet(t) = \prod_{\mathcal{I} \subset \{1, \ldots, n\}: \#\mathcal{I} = l-1} \zeta_{\mathcal{I} \cup \{0\}}^\bullet(t),$$

$$\zeta_\Lambda^\bullet(t) = \prod_{l=1}^{n+1} (\zeta_l^\bullet(t))^{(-1)^{l-1}}.$$

( $\bullet = 0$ or $\infty$) where $V_{l-1}(\gamma_s^\mathcal{I})$ is the (usual) $(l-1)$-dimensional volume of the face $\gamma_s^\mathcal{I}$ (in the hyperplane spanned by it in $L_{\mathcal{I} \cup \{0\}}$).

# 6 Examples

**Example 1.** Let $f = \frac{x^3 - xy}{y}$. The Milnor fibre $\mathcal{M}_f^0$ (repectively $\mathcal{M}_f^\infty$) is $\{(x,y) : \|(x,y)\| < \varepsilon, (x^3 - xy, y) \in B_\delta, x^3 - xy = cy\} \setminus \{(0,0)\}$, where $\|c\|$ is small (repectively large). From the equation $x^3 - xy = cy$ one has $y = \frac{x^3}{x+c}$ and thus $\mathcal{M}_f^0$ is diffeomorphic to the disk $\mathcal{D}$ in the $x$-plane with two points removed: $-c$ and the origin. In the same way $\mathcal{M}_f^\infty$ is diffeomorphic to the punctured disk $\mathcal{D}^*$. It is not difficult to understand that the action of the monodromy transformation in the homology groups is trivial in both cases. Thus

$$\zeta_f^0(t) = (1-t)^{-1} \quad \text{and} \quad \zeta_f^\infty(t) = 1.$$

Now let us calculate these zeta functions from their Newton diagrams, Fig 1.

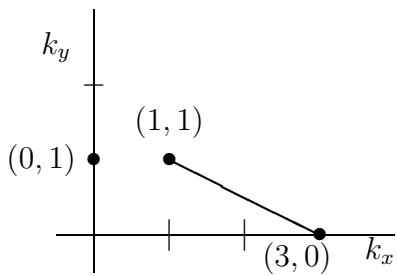

Figure 1.

We have $\zeta_1^\bullet(t) = 1$ since each coordinate axis intersects only one Newton diagram. There is only one linear function (namely $a = k_x + 2\,k_y$) such that $\dim \Delta(a, \Gamma_1) = 1$. The one-dimesional volume $V_1(\Delta(a, \Gamma_1))$ of $\Delta(a, \Gamma_1)$ is equal to 1 and $V_1(\Delta(a, \Gamma_2)) = 0$. We have $m(a, \Gamma_1) = 3$ and $m(a, \Gamma_2) = 2$. Thus $\zeta_2^0(t) = (1-t)$, $\zeta_2^\infty(t) = 1$, $\zeta_{(\Gamma_1,\Gamma_2)}^0(t) = (1-t)^{-1}$ and $\zeta_{(\Gamma_1,\Gamma_2)}^\infty(t) = 1$ which coincides with the formulae for $f$ written above.

**Example 2.** Let $P = xyz + x^p + y^q + z^r$ be a $T_{p,q,r}$ singularity, $\frac{1}{p} + \frac{1}{q} + \frac{1}{r} < 1$ and let $Q = x^d + y^d + z^d$ be a homogeneous polynomial of degree $d$. Suppose that $p > q > r > d > 3$ and that $p$, $q$, and $r$ are pairwise prime. Let us compute the zeta-functions of $f = \frac{P}{Q}$ using Theorems 2 and 3.

(a) It is clear that $f$ is non-degenerated with respect to its Newton pair $\Lambda = (\Gamma_1, \Gamma_2)$. Thus

$$\zeta_f^\bullet(t) = \zeta_\Lambda^\bullet(t) = \zeta_1^\bullet(\zeta_2^\bullet)^{-1}\zeta_3^\bullet \quad (\bullet = 0 \text{ or } \infty).$$

One has $\zeta_1^\infty = \zeta_2^\infty = 1$ and the unique covector which is necessary for computing $\zeta_3^\infty$ is $a = (1, 1, 1)$. In this case $m(a, \Gamma_1) = 3$, $m(a, \Gamma_2) = d$, $\Delta(a, \Gamma_1) = \{(1, 1, 1)\}$ and $\Delta(a, \Gamma_2)$ is the simplex $\{k_x + k_y + k_z = d, k_x \geq 0, k_y \geq 0, k_z \geq 0\}$, its two-dimensional volume is equal to $\frac{d^2}{2}$. Thus $\zeta_f^\infty = (1 - t^{d-3})^{d^2}$.

We have
$$\zeta_1^0 = (1 - t^{p-d})(1 - t^{q-d})(1 - t^{r-d}),$$
$$\zeta_2^0 = (1 - t^{r(q-d)})(1 - t^{r(p-d)})(1 - t^{q(p-d)})(1 - t^{r-d})^{2d}(1 - t^{q-d})^d.$$

To compute $\zeta_3^0$ one has to take into account both covectors $(rq - q - r, r, q)$, $(r, pr - p - r, p)$, and $(q, p, qp - p - q)$, corresponding to two-dimensional faces of $\Gamma_1$, and covectors $(1, r - 2, 1)$, $(r - 2, 1, 1)$, and $(q - 2, 1, 1)$, corresponding to pairs of the form (one-dimensional face of $\Gamma_1$, one-dimensional face of $\Gamma_2$). E.g., for $a = (1, r - 2, 1)$, $\Delta(a, \Gamma_1)$ (respectively $\Delta(a, \Gamma_2)$) is the segment betweenm $(0, 0, r)$ to $(1, 1, 1)$ (respectively between $(d, 0, 0)$ and $(0, 0, d)$. Pay attention to the "absence of the symmetry": last three covectors are not obtained from each other by permutting the coordinates and the numbers $p$, $q$, and $r$. This way

$$\zeta_3^0 = (1 - t^{r(q-d)})(1 - t^{r(p-d)})(1 - t^{q(p-d)})(1 - t^{r-d})^{2d}(1 - t^{q-d})^d$$

and

$$\zeta_f^0 = (1 - t^{p-d})(1 - t^{q-d})(1 - t^{r-d}).$$

(b) For computing the zeta-functions of $f$ with the help of Theorem 2, let $\pi : (\mathcal{X}, \mathcal{D}) \to (\mathbb{C}^3, 0)$ be the blowing-up of the origin in $\mathbb{C}^3$ and let $\varphi$ be the lifting $f \circ \pi$ of $f$ to the space $\mathcal{X}$. The exceptional divisor $\mathcal{D}$ is the complex projective plane $\mathbb{CP}^2$. Let $H_1$ and $H_2$ be the strict transforms of the hypersurfaces $\{P = 0\}$ and $\{Q = 0\}$, $D_i = \mathcal{D} \cap H_i$. The curve $D_1$ consists of three transversal lines $l_1, l_2, l_3$ and has three singular points $S_1 = l_2 \cap l_3 = (0, 0, 1)$, $S_2 = l_1 \cap l_3 = (0, 1, 0)$, and $S_3 = l_1 \cap l_3 = (1, 0, 0)$. The curve $D_2$ is a smooth curve of degree $d$, it intersects $D_1$ at $3d$ different points $\{P_1, \ldots, P_{3d}\}$.

One has the following natural stratification of the exceptional divisor $\mathcal{D}$:

(i) 0-dimensional strata $\Lambda_i^0$ ($i = 1, 2, 3$), each consisting of one point $S_i$;

(ii) 0-dimensional strata $\Xi_i^0$ consisting of one point $P_i$ each ($i = 1, \ldots, 3d$);

(iii) 1-dimensional strata $\Xi_i^1 = l_i \setminus \{D_2 \cup l_j \cup l_k\}$ ($i = 1, 2, 3$) and $\Xi_4^1 = D_2 \setminus D_1$;

(iv) 2-dimensional stratum $\Xi^2 = \mathcal{D} \setminus (D_1 \cup D_2)$.

It is not difficult to see that $\zeta_{\Xi^2}^0(t) = 1$, $\zeta_{\Xi^2}^\infty(t) = 1 - t^{d-3}$, for each stratum $\Xi$ from $\Xi_i^0$ ($1 \le i \le 3d$), $\Xi_i^1$ ($1 \le i \le 4$) one has $\zeta_\Xi^\bullet(t) = 1$ ($\bullet = 0$ or $\infty$).

In what follows the exceptional divisor $\mathcal{D}$ has the local equation $u = 0$. At the point $S_1$ the lifting $\varphi$ of the function $f$ is of the form $\frac{u^3 x_1 y_1 + u^r + x_1^p u^p + y_1^q u^q}{u^d x_1^d + u^d y_1^d + u^d}$. This germ has the same Newton pair as the germ $\frac{u^3 x_1 y_1 + u^r}{u^d}$. Using theorem 3 one has $\zeta_{\Lambda_1^0}^\infty = 1$, $\zeta_{\Lambda_1^0}^\infty = 1 - t^{r-d}$. At the point $S_2$ the function $\varphi$ has the form $\frac{u^3 x_1 z_1 + z_1^r u^r + x_1^p u^p + u^q}{u^d x_1^d + u^d + z_1^d u^d}$. It has the same Newton pair as $\frac{u^3 x_1 z_1 + z_1^r u^r + u^q}{u^d}$. Using Theorem 3 one has $\zeta_{\Lambda_2^0}^\infty(t) = 1$, $\zeta_{\Lambda_2^0}^0(t) = 1 - t^{q-d}$. Just in the same way $\zeta_{\Lambda_3^0}^\infty(t) = 1$, $\zeta_{\Lambda_3^0}^0(t) = 1 - t^{p-d}$. Combining these computations together, one has the same results as above (without using a partial resolution).

Addresses:
   Sabir M. Gusein-Zade,
Moscow State University,
Faculty of Mathematics and Mechanics.
Moscow, 119899, Russia.
   e-mail: sabir@ium.ips.ras.ru

   Ignacio Luengo,
Departamento de Algebra,
Facultad de Ciencias Matemáticas,
Universidad Complutense de Madrid.
Ciudad Universitaria, 28040 Madrid, Spain.
   e-mail: iluengo@eucmos.sim.ucm.es

   Alejandro Melle-Hernández,
Departamento de Geometría y Topología,
Facultad de Ciencias Matemáticas,
Universidad Complutense de Madrid.
Ciudad Universitaria, 28040 Madrid, Spain.
   e-mail: amelle@eucmos.sim.ucm.es